\documentclass[letterpaper, 10 pt]{ieeeconf}
\usepackage{cite}
\usepackage{amsmath,amssymb,amsfonts}
\usepackage{algorithmic}
\usepackage{graphicx}
\usepackage{textcomp}

\usepackage{tikz}
\usetikzlibrary{shapes,arrows}
\usepackage{nicefrac}
\usepackage{mathtools, cuted}
\usepackage{lipsum}
\usepackage{float}
\usepackage{afterpage}
\usepackage{placeins}
\usepackage{tabularx}
\usepackage{balance}
\usepackage{textcomp}%
\usepackage{booktabs}
\usepackage{multirow}
\usepackage{accents}
\usetikzlibrary{patterns,shapes,arrows}
\usetikzlibrary{positioning}
\allowdisplaybreaks

\newcommand\copyrighttext{%
  \footnotesize \textcopyright 2019 IEEE. Personal use of this material is permitted. Permission from IEEE must be obtained for all other uses, in any current or future media, including reprinting/republishing this material for advertising or promotional purposes, creating new collective works, for resale or redistribution to servers or lists, or reuse of any copyrighted component of this work in other works.}
\newcommand\copyrightnotice{%
\begin{tikzpicture}[remember picture,overlay]
\node[anchor=south,yshift=10pt] at (current page.south) {\fbox{\parbox{\dimexpr\textwidth-\fboxsep-\fboxrule\relax}{\copyrighttext}}};
\end{tikzpicture}%
}

\tikzstyle{decision} = [diamond, draw, fill=blue!20,
    text width=7em, text badly centered, inner sep=0pt]
\tikzstyle{decision1} = [diamond, draw, fill=blue!20,
    text width=7em, text badly centered, inner sep=0pt, node distance = 2.5cm]
\tikzstyle{block} = [rectangle, draw, fill=blue!20,
    text width=13em, text centered, rounded corners, minimum height=3em]
\tikzstyle{block2} = [rectangle, draw, fill=blue!20,
    text width=5em, text centered, rounded corners, minimum height=3em, node distance = 4cm]
\tikzstyle{block3} = [rectangle, draw, fill=red!20,
    text width=5em, text centered, rounded corners, minimum height=3em, node distance = 4cm]
\tikzstyle{line} = [draw, very thick, color=black!50, -latex']
\tikzstyle{cloud} = [draw, ellipse,fill=red!20, node distance=2.5cm, minimum height=2em]

\newtheorem{remark}{Remark}

\begin{document}

\title{Efficient and More Accurate Representation of Solution Trajectories in Numerical Optimal Control}

\author{Yuanbo Nie and Eric C.\ Kerrigan \IEEEmembership{Senior Member, IEEE}
\thanks{Yuanbo Nie and Eric C.\ Kerrigan are with the Department of Aeronautics, Imperial College London, SW7~2AZ, U.K. {\tt\small yn15@ic.ac.uk}, {\tt\small 
e.kerrigan@imperial.ac.uk}}%
\thanks{Eric C. Kerrigan is also with the Department of Electrical \& Electronic Engineering, Imperial College London, London SW7~2AZ, U.K.}%
\thanks{Accepted version to be published in: IEEE Control Systems Letters}%
}

\maketitle
\copyrightnotice

\begin{abstract}
We show via examples that, when solving optimal control problems, representing the optimal state and input trajectory directly using interpolation schemes may not be the best choice. Due to the lack of considerations for solution trajectories  in-between  collocation points, large errors may occur, posing risks if this solution is to be  applied. A novel solution representation method is proposed, capable of yielding a solution of much higher accuracy for the same discretization mesh. This is achieved by minimizing the integral of the residual error for the overall trajectory, instead of forcing the errors to be zero only at collocation points. In this way, the requirement for mesh resolution can be significantly reduced, leaving the problem dimensions relatively small.  This particular formulation also avoids some of the drawbacks found in the earlier work of integrated residual minimization, leading to more efficient computations. 
\end{abstract}

\begin{keywords}
residual minimization, optimal control, solution representation, nonlinear predictive control
\end{keywords}

\section{Introduction} \label{sec:Intro}
Solving optimal control problems (OCPs) are central in the field of trajectory optimization and real-time optimization-based control.  In practice, most OCPs need to be solved with numerical schemes. In many situations, indirect methods can be difficult to implement, since they require analytic expressions of optimality conditions, which can be hard to derive. Direct methods have consequently become the standard for solving practical optimal control problems~\cite{limebeer2015faster}. 

Direct methods often require transcription of the infinite-dimensional OCP into a nonlinear programming (NLP) problem of a finite dimension, via the introduction of a discretization mesh. If only the control trajectories are discretized as decision variables, the method is called sequential, of which direct single shooting is an example. In simultaneous approaches, such as direct multiple shooting and direct collocation, the state trajectories are also discretized and included as decision variables. A comprehensive comparison of different solution strategies are available in \cite{cervantes2001optimization,biegler2007overview}, demonstrating the advantages of simultaneous methods in solving large-scale problems using sparse NLP formulations, compared to sequential methods.


Since the NLP solver will only return samples of the solution at a finite number of points in the discretization mesh, additional steps and care must be taken when representing the continuous-time results for time instances other than these sampled points. The common practice today is to directly use interpolation schemes, which are selected in accordance with the type of discretization mesh \cite{betts2010practical,kelly2017introduction}. 

Collocation methods with direct interpolation have a major drawback: the  residual error of the ordinary differential equations (ODEs)  is forced to zero only at collocation points. In general, no guarantees on accuracy and constraint satisfaction can be derived for system trajectories  in-between  collocation points. Posterior analysis is needed to identify intervals where errors are high and mesh refinement procedures are used to modify the discretization mesh. The problem has to be solved iteratively until all the errors are within tolerances.

We present a solution representation method that can significantly improve the solution accuracy for the same discretization mesh, so that results of higher quality are obtainable with relatively coarse meshes. This is achieved by minimizing the ODE residual error integrated over the solution trajectory. Section~\ref{sec: OptimizationBasedControl} provides the background information for solving optimal control problems numerically with the direct collocation method. Section~\ref{sec: ResidualMinimization} introduces the fundamental concept of the residual minimization method and motivates the development of the proposed scheme, which is presented in Section~\ref{sec: ProposedScheme}. The benefits of the method are demonstrated in Section~\ref{sec: ExampleProblem} with two example problems, followed by concluding remarks in Section~\ref{sec: conlclusions}.

\section{Numerical Optimal Control}
\label{sec: OptimizationBasedControl}

Optimization-based control often requires the solution of  OCPs expressed in the general Bolza form:
\begin{subequations}
\label{eqn: OCPBolza}
\begin{equation}
\min_{x,u,p,t_0,t_f} \Phi(x(t_0),t_0,x(t_f),t_f,p)
+\int_{t_0}^{t_f} L(x(t),u(t),t,p)\: dt
\end{equation}
subject to
\begin{align}
\dot{x}(t)=f(x(t),u(t),t,p),\ &\forall t \in [t_0,t_f] \label{eqn:OCPBolzaDynamics}\\
c(x(t),u(t),t,p)\le 0,\ &\forall t \in [t_0,t_f] \label{eqn:OCPBolzaPathConstraint}\\
\phi(x(t_0),t_0,x(t_f),t_f,p) =0,\ &
\end{align}
\end{subequations}
 with $x: \mathbb{R} \rightarrow \mathbb{R}^n$ is the state trajectory of the system, $u: \mathbb{R} \rightarrow \mathbb{R}^m$ is the control input trajectory,   $p \in \mathbb{R}^s$ are static parameters, $t_0 \in \mathbb{R}$ and $t_f \in \mathbb{R}$ are the initial and terminal time.  $\Phi$ is the Mayer cost functional ($\Phi$: $\mathbb{R}^n \times \mathbb{R} \times \mathbb{R}^n \times \mathbb{R} \times \mathbb{R}^s \to \mathbb{R}$), $L$ is the Lagrange cost functional ($L:\mathbb{R}^n \times \mathbb{R}^m \times \mathbb{R} \times \mathbb{R}^s \to \mathbb{R}$), $f$ is the dynamic constraint ($f:\mathbb{R}^n \times \mathbb{R}^m \times \mathbb{R} \times \mathbb{R}^s \to \mathbb{R}^n$), $c$ is the path constraint ($c:\mathbb{R}^n \times \mathbb{R}^m \times \mathbb{R} \times \mathbb{R}^s \to \mathbb{R}^{n_g}$) and $\phi$ is the boundary condition ($\phi:\mathbb{R}^n \times \mathbb{R} \times \mathbb{R}^n \times \mathbb{R} \times \mathbb{R}^s \to \mathbb{R}^{n_q}$).

\subsection{Direct collocation methods}
\label{sec: DirectTranscriptionMethod}
Direct collocation methods can be categorized into fixed-order $h$ methods (e.g.\ Euler, Trapezoidal, and Hermite-Simpson (H-S) as in \cite{betts2010practical}), and variable higher-order $p$/$hp$ methods (e.g.\ Legendre-Gauss-Radau (LGR) as in \cite{liu2014hp}). Here, we aim to provide a high level overview. With a mesh of size $N:=\sum_{k=1}^K N^{(k)}$, the states can be approximated as

\begin{equation}
\label{eqn: LGRStateApproximation}
x^{(k)}(\tau) \approx \bar{x}^{(k)}(\tau) := \sum_{j=1}^{N^{(k)}}\mathcal{X}_j^{(k)}\mathcal{B}_{j}^{(k)}(\tau),
\end{equation}

within mesh interval $k$ $\in$ $\{1,\ldots, K\}$, where $N^{(k)}$ is the number of collocation points  for  interval $k$, and $\mathcal{B}_{j}^{(k)}(\cdot)$ are basis functions. For classical $h$ methods, $\tau \in \mathbb{R}^{N}$ takes on values in the interval $[0,1]$ representing $[t_0,t_f]$, and $\mathcal{B}_{j}^{(k)}(\cdot)$ are chosen to be elementary B-splines of various orders. For $p$/$hp$ methods, $\mathcal{B}_{j}^{(k)}(\cdot)$ are Lagrange interpolating polynomials over the normalized time interval $\tau$ $\in$ $[-1,1]$. We use $X_j^{(k)}$ and $U_j^{(k)}$ to represent the approximated states and inputs at collocation points, e.g.\  $X_j^{(k)}=\bar{x}^{(k)}(\tau_j^{(k)}) \in \mathbb{R}^{n}$, where $\tau_j^{(k)}$ is the $j^\text{th}$ collocation point in mesh interval~$k$.

Consequently, the OCP~\eqref{eqn: OCPBolza} can be approximated by
\begin{subequations}
\label{eqn: LGRStateApproximationAll}
\begin{multline}
\label{eqn: LGRStateApproximationCost}
J_c:=\min_{X,U,p,t_0,t_f}  \Phi(X_1^{(1)},t_0,X_{f}^{(K)},t_f,p)\\
+\sum_{k=1}^{K}\sum_{i=1}^{N^{(k)}} w_i^{(k)} L(X_i^{(k)},U_i^{(k)},\tau_i^{(k)},t_0,t_f,p)
\end{multline}
subject to, for $i=1,\hdots,N^{(k)}$ and $k=1,\hdots,K$:
\begin{align}
\label{eqn: LGRStateApproximationCostDefect}
\sum_{j=1}^{N^{(k)}}\mathcal{A}_{ij}^{(k)}X_j^{(k)}+\mathcal{D}_{i}^{(k)}f(X_i^{(k)},U_i^{(k)},\tau_i^{(k)},t_0,t_f,p) = & 0 \\
\label{eqn: LGRStateApproximationPathConstraint}
c(X_i^{(k)},U_i^{(k)},\tau_i^{(k)},t_0,t_f,p)\le & 0  \\
\phi(X_1^{(1)},t_0,X_{f}^{(K)},t_f,p) =& 0
\end{align}
\end{subequations}
where $w_j^{(k)}$ are the quadrature weights for the respective discretization method chosen, $\mathcal{A}$ is the numerical differentiation matrix with $\mathcal{A}_{ij}$ the element $(i,j)$ of the matrix, and $\mathcal{D}$ a constant matrix. The discretized problem can then be solved with off-the-shelf or structure-exploiting NLP solvers.

\subsection{Representing the results}
\label{subsec: ResultReconstruction}
The NLP solver generates a discretized solution $\mathcal{Z} \coloneqq (X, U, p, t_0, t_f)$ as sampled data points.  Interpolating splines may be used to construct an approximation of the continuous-time solution $\tilde{z}(t) \coloneqq (\tilde{x}(\mathcal{Z},t), \tilde{u}(\mathcal{Z},t), t, p)$, with $\tilde{x}(\mathcal{Z},\cdot)$, $\tilde{u}(\mathcal{Z},\cdot)$ the approximated state and input trajectories. 
%
%
%

\subsubsection{Representation via direct interpolation}
\label{subsec: ReconstructionDefault}
Conventionally, the interpolation of the solution corresponds to the discretization scheme used in the transcription process. Thus, we must analyze how the state approximation \eqref{eqn: LGRStateApproximation} enters the optimal control problem formulation~\eqref{eqn: LGRStateApproximationCost}.

It is not difficult to discover that the only dependency on the basis function in \eqref{eqn: LGRStateApproximationCost} appears in the defect constraint~\eqref{eqn: LGRStateApproximationCostDefect}, through the first term  representing the numerical differentiation of the approximated function  $\bar{x}(\cdot)$.  

For most commonly-used numerical schemes, the numerical differentiation formulation has an equivalent integration form. Both forms are presented in Table \ref{tab: NumericalIntegrationSchemes}, where $h_k := \Delta t(\tau_{N}^{(k)}-\tau_1^{(k)})$, $\Delta t := t_f-t_0$, and 
\begin{equation}
\label{eqn: Fi_Discretized}
F_i^{(k)}:=f(X_i^{(k)},U_i^{(k)},\tau_i^{(k)},t_0,t_f,p).
\end{equation}

\begin{table}[b]
  \small
	\begin{center}
	\caption{Typical numerical schemes}
	\label{tab: NumericalIntegrationSchemes}
		\begin{tabular}{c|c}
		 \textbf{Method}   & \multirow{2}{*}{\textbf{Numerical Integration Scheme}} \\
		\textbf{(Order)} & \\
		\hline \multirow{2}{*}{Euler (1)} & \multirow{2}{*}{$X_{2}^{(k)}=X_1^{(k)}+h_kF_1^{(k)}$}\\
		& \\
		\hline \multirow{2}{*}{Trapezoidal (2)} & \multirow{2}{*}{$X_{2}^{(k)}=X_1^{(k)}+\frac{h_k}{2}(F_1^{(k)}+F_{2}^{(k)})$} \\
		 & \\
		\hline Hermite & $X_{2}^{(k)}=\frac{1}{2}(X_{2}^{(k)}+X_1^{(k)})+\frac{h_k}{8}(F_1^{(k)}-F_{3}^{(k)})$\\
		 Simpson (3)& $X_{3}^{(k)}=X_1^{(k)}+\frac{h_k}{6}(F_1^{(k)}+4F_{2}^{(k)}+F_{3}^{(k)})$\\
		\hline \multirow{2}{*}{LGR ($N^{(k)}$)} & $\mathcal{I}^{(k)}=[\mathcal{A}^{(k)}_{2:N+1}]^{-1}$ \\
		 & $X_{2:N+1}^{(k)}=X_1+\frac{\Delta t}{2}\mathcal{I}^{(k)}F_{1:N}^{(k)}$\\
		\hline
		\end{tabular} 
	\end{center}
\end{table}

For each numerical scheme, direct interpolation of the OCP solution is possible using splines with the type and order in accordance with Table \ref{tab: ReconstructionContinuity}. 
\begin{table}[tb]
  \small
	\begin{center}
	\caption{Continuity of the reconstructed solution (a.m.: at most; p.w.: piecewise)}
	\label{tab: ReconstructionContinuity}
		\begin{tabular}{c|c|c|c}
		\hline \textbf{Method} & \textbf{Dynamics ($\dot{\tilde{x}}$)} & \textbf{States ($\tilde{x}$)}  & \textbf{Inputs ($\tilde{u}$)} \\
		\hline Euler & p.w.\ constant & a.m.p.w.\ linear & \multirow{4}{*}{\shortstack{same \\ as \\dynamics}} \\
		\cline{1-3} Trape. & a.m.p.w.\ linear & a.m.p.w.\ quad. &\\
		\cline{1-3} H-S & a.m.p.w.\ quad. & a.m.p.w.\ cubic & \\
		\cline{1-3} LGR & a.m.\  order $N^{(k)}$ & a.m.\ order $N^{(k)}$+1 &\\
		\hline
		\end{tabular} 
	\end{center}
\end{table}
For example, with Hermite-Simpson transcription, the reconstructed state trajectory inside mesh interval $k$ using cubic splines will be
\begin{multline}
\label{eqn: HSStateReconstruction}
\tilde{x}^{(k)}(\mathcal{Z},t)=X_1^{(k)}+F_1^{(k)}(t-t_1^{(k)})\\
+\frac{1}{2}\bigg(-3F_1^{(k)}+4F_2^{(k)}-F_3^{(k)}\bigg)\frac{(t-t_1^{(k)})^2}{h_k}\\
+\frac{2}{3}\bigg(F_1^{(k)}-2F_2^{(k)}+F_3^{(k)}\bigg)\frac{(t-t_1^{(k)})^3}{h_k^2}, 
\end{multline}
the dynamics trajectory with quadratic splines will be
\begin{multline}
\label{eqn: HSDynamicsReconstruction}
\dot{\tilde{x}}^{(k)}(\mathcal{Z},t)=F_1^{(k)}+\bigg(-3F_1^{(k)}+4F_2^{(k)}-F_3^{(k)}\bigg)\frac{t-t_1^{(k)}}{h_k}\\
+\bigg(2F_1^{(k)}-4F_2^{(k)}+2F_3^{(k)}\bigg)\bigg(\frac{t-t_1^{(k)}}{h_k}\bigg)^2, 
\end{multline}
and the control trajectory with quadratic splines will have the expression
\begin{multline}
\label{eqn: HSInputReconstruction}
\tilde{u}^{(k)}(\mathcal{Z},t)=\frac{2}{h_k^2}(t-\frac{1}{2}t_1^{(k)}-\frac{1}{2}t_3^{(k)})(t-t_3^{(k)})U_1^{(k)}\\
-\frac{4}{h_k^2}(t-t_1^{(k)})(t-t_3^{(k)})U_2^{(k)}\\
+\frac{2}{h_k^2}(t-t_1^{(k)})(t-\frac{1}{2}t_1^{(k)}-\frac{1}{2}t_3^{(k)})U_3^{(k)}.
\end{multline}
for all $t \in [t_1^{(k)}, t_3^{(k)}]$. The whole trajectory $\tilde{x}(\mathcal{Z},\cdot)$, $\dot{\tilde{x}}(\mathcal{Z},\cdot)$ and $\tilde{u}(\mathcal{Z},\cdot)$ can then be expressed as piecewise polynomials.

For $p/hp$ methods, Lagrange interpolating polynomials are often used as basis functions during the transcription process. An alternative version, namely barycentric Lagrange interpolation, is often used instead for solution interpolation, due to its improved numerical stability.

\subsubsection{Evaluation of errors}
\label{subsec: SolRepresentationAndErrorAnalysis}

The quality of the interpolated solution needs to be assured through error analysis, assessing the level of accuracy and constraint satisfaction. Firstly, any valid trajectory $\tilde{z}(\cdot)$ must satisfy the system dynamics \eqref{eqn:OCPBolzaDynamics} with a good level of accuracy. Therefore, one measure for the error due to discretization and interpolation is through the calculation of the ODE residual $\varepsilon_r(t) \in \mathbb{R}^n$ defined as
\begin{equation}
\label{eqn: DiscretizationError}
\varepsilon_r(t):=\dot{\tilde{x}}(\mathcal{Z},t)-f(\tilde{x}(\mathcal{Z},t), \tilde{u}(\mathcal{Z},t), t, p).
\end{equation}
For the discretized problem, the error in the state variables over each interval  in-between  collocation points can then be estimated with the integral
\begin{equation*}
\eta_j:=\int^{t_{j+1}}_{t_j} 	\|\varepsilon_r(s)\|_{2}\: ds,
\end{equation*}
as a single metric for a multi-variable problem, or
\begin{equation*}
\sigma_{j,q}:=\int^{t_{j+1}}_{t_j} |\varepsilon_{r_q}(s)|\: ds, \text{ for } q=1,\hdots,n,
\end{equation*}
for each dynamics equation separately. $\eta \in \mathbb{R}^N$ or $\sigma \in \mathbb{R}^{n \times N}$ are typically referred to as the \emph{absolute local error}~\cite{betts2010practical}. The operator $\|\cdot\|_{2}$ is the vector 2-norm. The integral can be practically estimated by high order quadrature. 

In addition, numerical discretization inevitably leads to possible constraint violations of the trajectories  in-between  the collocation points. For path and box constraints that are expressed semi-explicitly as \eqref{eqn:OCPBolzaPathConstraint}, the \textit{absolute local constraint violation} $\varepsilon_{c_\zeta}(t) \in \mathbb{R}^{n_g}$ may be straight-forwardly estimated by
\begin{align*}
\label{eqn: ConstraintViolationError}
\varepsilon_{c_\zeta}(t):=\begin{cases}
0 & \text{if } c_\zeta(\tilde{z}(t)) \leq 0\\ 
c_\zeta(\tilde{z}(t)) & \text{if } c_\zeta(\tilde{z}(t)) > 0\\
\end{cases}, \text{ for } \zeta=1,\hdots,n_g.
\end{align*}

Once the distributions of errors are calculated, appropriate modifications can be made to the discretization mesh, to iteratively resolve the problem until the obtained solution fulfills all predefined error tolerances ($\eta_{tol}$ and $\varepsilon_{c_{tol}}$). This process is called mesh refinement (MR). Common approaches for mesh refinement include adding intervals and/or changing the polynomial order. The NLP formulated based on the new mesh is warm started using the previous solution from the coarser mesh. This can often lead to significantly faster convergence, compared to a fine uniform mesh without MR, thus reducing the overall computation time.

\subsection{Problems associated with direct reconstruction}
Practical experience has shown that trajectory interpolation in accordance with the discretization scheme is not the best choice. In many cases  large discretization errors and constraint violations occur inside the intervals  in-between  collocation points. Furthermore, if the optimal control trajectory is discontinuous, direct interpolation using polynomials can often result in a Gibbs-like phenomenon, inducing non-physical oscillations in the solution.  A typical example for this to happen would be in problems with bang-bang control, if the switch happens  in-between mesh intervals. 

These issues are fundamentally rooted in the direct collocation formulation. Firstly, states, dynamics and controls can rarely all be approximated accurately by polynomials. Even in the simple case where $f(x(t),u(t))=\dot{x}(t)=ax(t)+u(t)$ and $u(t)=1$ are both polynomials (thus can be represented exactly by polynomials), the corresponding state trajectory $x(t)=x(0)e^{at}+\int_0^t e^{a(t-s)}u(s)\ ds$ is clearly not a polynomial and approximation errors should be expected. 

It is then important to note that driving the defect constraint \eqref{eqn: LGRStateApproximationCostDefect} to zero (or machine precision) at collocation points does not imply that the polynomial functions used for the state and input approximations in the NLP will  satisfy the dynamic equations and constraints   in-between  collocation points. In fact, the opposite can and often does occur.

It is well-known in the field of curve fitting that if a function cannot be exactly represented by a polynomial, forcing the polynomial to exactly go though some sampled data points generally results in larger errors in comparison to fitting using least squares criteria. The same analogy can be applied here: forcing the defect constraints to be zero at collocation points will generally result in larger overall defect errors for the whole trajectory, in comparison to a method that minimizes the integral of the defect errors in a least squares manner. This observation motivated the development of the integrated residual minimization method.

\section{Method of integrated residual minimization}
\label{sec: ResidualMinimization}
Integrated residual minimization is motivated by the recently-proposed method in~\cite{neuenhofen2018dynamic}, which is a generalization of the least-squares approach for solving differential equations to solving dynamic optimization problems. The idea is that instead of forcing the ODE residuals~\eqref{eqn: DiscretizationError} to be zero at collocation points with \eqref{eqn: LGRStateApproximationCostDefect}, the method tries to minimize the square of the 2-norm of the ODE residuals for the represented solution polynomials integrated along the whole trajectory, i.e.
\begin{equation}
\label{eqn: ResidualMinimizationOrg}
\min_{\hat{x}, \hat{u}, t, p} \int_{t_0}^{t_f} r(\hat{x}(t), \hat{u}(t), t, p)\: dt
\end{equation}
with
\begin{equation}
\label{eqn: ResidualMinimizationOrgr}
r(\hat{x}(t), \hat{u}(t), t, p):=\| \dot{\hat{x}}(t)-f(\hat{x}(t), \hat{u}(t), t, p)\|^2_{2}.
\end{equation}
As presented in \cite{neuenhofen2018dynamic}, the expressions for state and input functions $\hat{x}$ and $\hat{u}$ can be polynomials of any standard types, with polynomial coefficients $P_{j,q}$ as decision variables. 

 This concept shares some resemblance with the method of direct error enforcement in \cite{vasantharajan1990simultaneous}, where the residual errors are computed with a different formulation and only at a single non-collocation reference point for each mesh interval. 

\begin{remark}
\label{rem: shortcomingRM}
The choice of representation in \cite{neuenhofen2018dynamic} with $P_{j,q}$ as decision variables increases the computational complexity of the problem in comparison to direct collocation: 
\begin{itemize}
\item One extra decision variable is required for every state and input variable in every mesh segment.
\item Simple bounds need to be implemented as general inequality constraints.
\item Additional computations are required to obtain the initial guesses of the decision variables from an estimation of the solution trajectory.
\item The magnitudes of decision variables may span a wide numerical range. This is detrimental in terms of ensuring consistent numerical accuracy in computations. 
\item If finite differences are used for obtaining the derivative information, the calculations can be less accurate.
\item Proper scaling of decision variables can be difficult.
\item State continuity in-between mesh segments might need to be enforced with additional equality constraints.
\end{itemize}
\end{remark}

Therefore, we need to develop a method that avoids the above-listed drawbacks and, to a great extent, retains the computational efficiency of direct collocation.

\section{The proposed scheme}
\label{sec: ProposedScheme}
Based on the above observations, we propose a method to generate solution trajectories that can be  orders of magnitude  more accurate than direct interpolation for the ODE defect error, without increasing the size of the discretization mesh. The method retains the same decision variables as in \eqref{eqn: LGRStateApproximationAll}, namely $\mathcal{Z} \coloneqq (X, U, p, t_0, t_f)$, and uses the interpolation polynomial formula $\tilde{x}(\mathcal{Z},\cdot)$, $\dot{\tilde{x}}(\mathcal{Z},\cdot)$ and $\tilde{u}(\mathcal{Z},\cdot)$ to directly  represent  $\hat{x}(\cdot)$, $\dot{\hat{x}}(\cdot)$ and $\hat{u}(\cdot)$ in \eqref{eqn: ResidualMinimizationOrg} and \eqref{eqn: ResidualMinimizationOrgr}. 

For example, consider Hermite-Simpson discretization. The input trajectory inside mesh interval $k$ can be represented by the polynomial as in \eqref{eqn: HSInputReconstruction}, based on the values of the decision variables $U_1^{(k)}$, $U_2^{(k)}$ and $U_3^{(k)}$. However, for the state trajectory, one challenge arises. For solutions to \eqref{eqn: LGRStateApproximationAll}, continuity of state variables are automatically fulfilled when using \eqref{eqn: HSStateReconstruction} as the interpolation equation; however, this is not generally the case for arbitrary solutions that violate the defect constraint~\eqref{eqn: LGRStateApproximationCostDefect}. 

To avoid imposing additional path constraints for state continuity, we make use of the original Hermite-Simpson numerical integration scheme (in Table \ref{tab: NumericalIntegrationSchemes}), and obtain 
\begin{align}
\label{eqn: HSIntegrationF2}
F_2^{(k)}=&-\frac{1}{2h_k}(5X_1^{(k)}-4X_2^{(k)}-X_3^{(k)}+F_1^{(k)}h_k)\\
\label{eqn: HSIntegrationF3}
F_3^{(k)}=&\frac{1}{h_k}(4X_1^{(k)}-8X_2^{(k)}+4X_3^{(k)}+F_1^{(k)}h_k).
\end{align}
As a check, substituting $t=t_2^{(k)}=t_1^{(k)}+h_k/2$ and \eqref{eqn: HSIntegrationF2} into~\eqref{eqn: HSStateReconstruction} will result in $X_2^{(k)}$, and substituting $t=t_3^{(k)}=t_1^{(k)}+h_k$ and \eqref{eqn: HSIntegrationF3} into \eqref{eqn: HSStateReconstruction} will result in $X_3^{(k)}$. Thus, with $F_1^{(k)}$, $F_2^{(k)}$ and $F_3^{(k)}$ calculated based on \eqref{eqn: Fi_Discretized}, \eqref{eqn: HSIntegrationF2} and \eqref{eqn: HSIntegrationF3}, respectively, the interpolation formula \eqref{eqn: HSStateReconstruction} guarantees state trajectory continuity without imposing additional constraints. 

Thus, an optimization problem for representing the OCP solution can be formulated as
\begin{subequations}
\label{eqn: ResMinInterpolationAll}
\begin{equation}
\label{eqn: ResMinInterpolationCost}
\min_{X,U,p,t_0,t_f} \sum_{k=1}^{K} R(X^{(k)},U^{(k)},\tau^{(k)},\tau_{q}^{(k)},t_0,t_f,p)
\end{equation}
subject to, for $i=1,\hdots,N^{(k)}$ and $k=1,\hdots,K$,
\begin{align}
\begin{split}
\label{eqn: ResMinInterpolationOrgCost}
\sum_{k=1}^{K}\sum_{i=1}^{N^{(k)}} w_i^{(k)} L(X_i^{(k)},U_i^{(k)},\tau_i^{(k)},t_0,t_f,p) \quad &\\
+\Phi(X_1^{(1)},t_0,X_{f}^{(K)},t_f,p) \le & J_c  
\end{split}\\
\label{eqn: ResMinInterpolationPathConstraint}
c(X_i^{(k)},U_i^{(k)},\tau_i^{(k)},t_0,t_f,p)\le & 0  \\
\phi(X_1^{(1)},t_0,X_{f}^{(K)},t_f,p) =& 0
\end{align}
\end{subequations}
with  \eqref{eqn: ResMinInterpolationOrgCost} the constraint for the objective of the original OCP~\eqref{eqn: LGRStateApproximationAll} and  $J_c \in \mathbb{R}$ the value of the cost obtained from direct collocation. $R$ is the residual cost: for certain problems, this can be calculated precisely with analytical expressions;  for most practical problems, quadrature rules of sufficiently high order can be used, i.e.
\begin{multline*}
R(X^{(k)},U^{(k)},\tau^{(k)},\tau_{q}^{(k)},t_0,t_f,p):=\\
\sum_{\iota=1}^{N_q^{(k)}}w_{\iota}^{(k)}\|\dot{\tilde{x}}(\mathcal{Z},t_{q_{\iota}}^{(k)})-f(\tilde{x}(\mathcal{Z},t_{q_{\iota}}^{(k)}), \tilde{u}(\mathcal{Z},t_{q_{\iota}}^{(k)}), t_{q_{\iota}}^{(k)}, p)\|^2_{2}
\end{multline*}
with $t_{q_{\iota}}^{(k)}:=\frac{t_f^{(k)}-t_0^{(k)}}{2}\tau_{q_{\iota}}^{(k)}+\frac{t_f^{(k)}+t_0^{(k)}}{2}$ and $\tau_q^{(k)} \in \mathbb{R}^{N_q^{(k)}}$ the quadrature mesh for approximating the integral inside a mesh interval, where $w_{\iota}^{(k)}$ are the corresponding quadrature weights. Typically a Gaussian quadrature of order $N_q^{(k)}\ge 4N^{(k)}+1$ is required for a good accuracy \cite{neuenhofen2018dynamic}. 

 When comparing \eqref{eqn: ResMinInterpolationAll} to the original OCP formulation of direct collocation \eqref{eqn: LGRStateApproximationAll}, it is straightforward to note that, if $\mathcal{Z}$ is a solution to \eqref{eqn: LGRStateApproximationAll}, then $\mathcal{Z}$ will also be a feasible point for~\eqref{eqn: ResMinInterpolationAll}. If no better feasible solution is obtainable, \eqref{eqn: ResMinInterpolationAll} is at least guaranteed to have one solution, namely $\mathcal{Z}$.

In terms of computational complexity, the proposed formulation \eqref{eqn: ResMinInterpolationAll} avoids the shortcomings as listed in Remark \ref{rem: shortcomingRM}.  In addition, unlike the penalty-barrier finite element (PBF) method proposed in \cite{neuenhofen2018dynamic}, which requires tailored solvers for a good performance, \eqref{eqn: ResMinInterpolationAll} can be efficiently solved with the same off-the-shelf sparse NLP solvers as direct collocation. The transcription and the majority of the computational components can be shared between the two, and warm starting techniques can be exploited to  accelerate the computations.

\section{Example Problems}
\label{sec: ExampleProblem}
Here, we present two example problems to demonstrate the main advantages of the proposed scheme. Both OCPs are transcribed using the optimal control software \texttt{ICLOCS2}~\cite{ICLOCS2}, and numerically solved to a tolerance of $10^{-9}$ with NLP solver \texttt{IPOPT} \cite{wachter2006implementation} (version 3.12.9). With extremely coarse meshes, the emphasis of the comparison will not be on yielding solutions that look similar to the optimal trajectory.  Instead, the goal is to obtain sub-optimal solutions that, when applied, can result in low discrepancies between the represented solution and the implementation outcome.

\subsection{Two-Link Robot Arm}

The two-link robot arm problem presented here was adapted from \cite[Ex.\ 2, Sect.\ 12.4.2]{luus2000iterative}. Consider a system consisting of two identical beams with the same property (mass: $m=1$\,kg, length: $l=1$\,m, and moment of inertia), connected at two actuated joints. The objective is to reposition  a payload of mass $M=1$\,kg in minimum time, with the addition of a regularization term:
\begin{equation*}
\min_{x,u,t_f} \quad t_f+0.01\int_{0}^{t_f} u_1(t)^2+u_2(t)^2\: dt. \\
\end{equation*}
The system has angular rates $\omega_{\phi}$, $\omega_{\psi}$, and angles $\phi$,  $\chi=\phi-\psi$  as state variables, and nondimensionalized torque $u_1$ and $u_2$ as inputs. Furthermore, the variable simple bounds and boundary conditions are imposed in accordance to the reference, except that $\chi(t_f)=0.5$\,rad, and $\phi(t_f)=0.522$\,rad.

\begin{figure}[t]
\begin{center}
\includegraphics[width=\columnwidth]{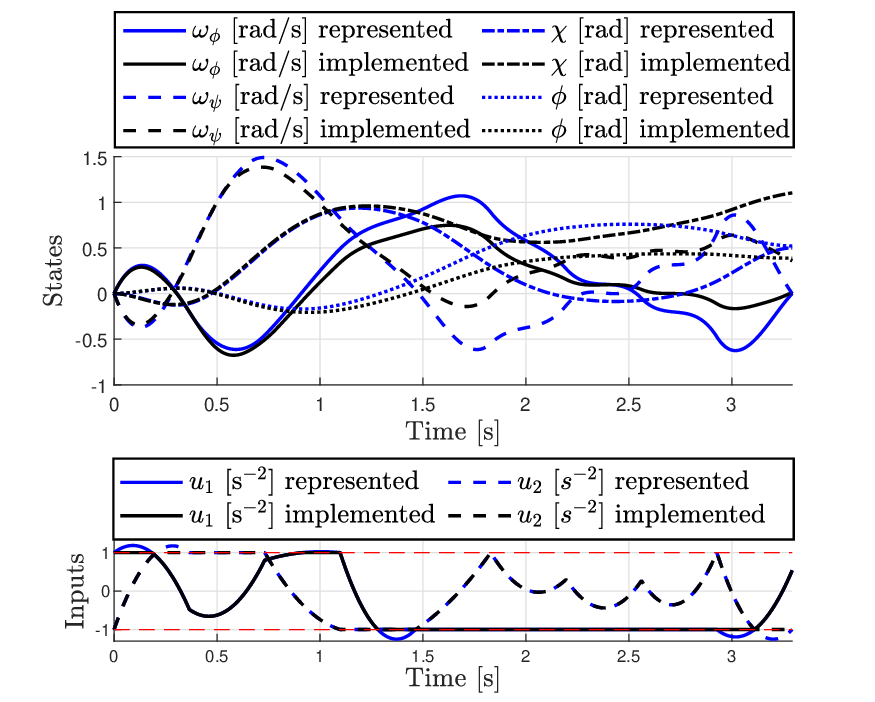}    
\caption{Solution to the two-link robot arm problem, direct interpolation method for solution representation, direct collocation with Hermite-Simpson discretization, 10 mesh intervals} 
\label{fig:Coll_roboticArm}
\end{center}
\end{figure}

\begin{figure}[t]
\begin{center}
\includegraphics[width=\columnwidth]{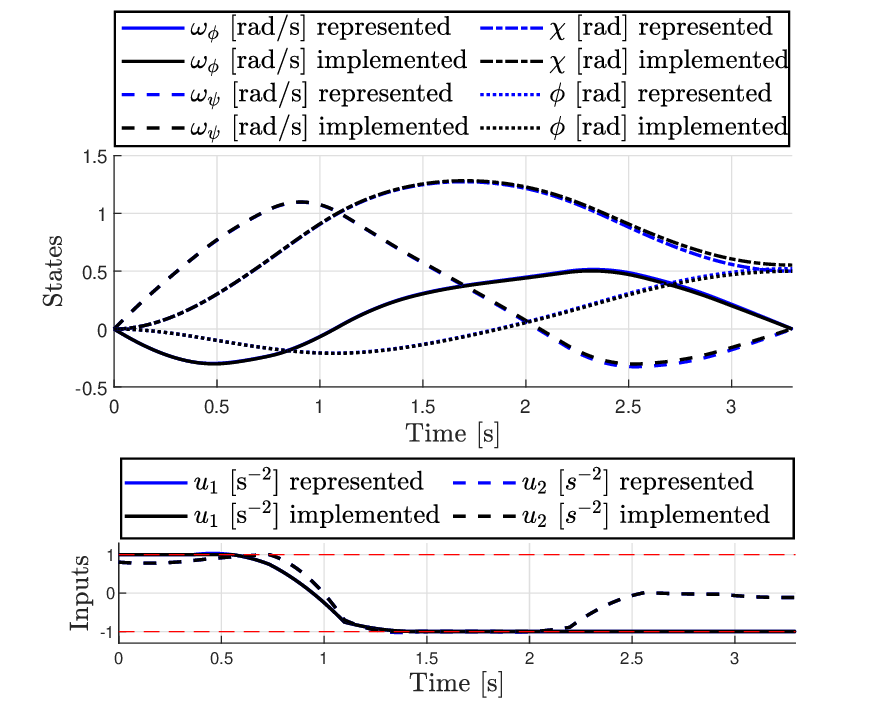}    
\caption{Solution to the two-link robot arm problem, integrated residual minimization method for solution representation, direct collocation with Hermite-Simpson discretization, 10 mesh intervals} 
\label{fig:ResMin_roboticArm}
\end{center}
\end{figure}

Figures \ref{fig:Coll_roboticArm} and \ref{fig:ResMin_roboticArm} illustrate the solutions to the two-link robot arm problem problem generated with the two different solution representation methods. 
Presented alongside are the outcomes from the actual implementation of the resultant input trajectory on the same dynamic model, solved with a non-stiff variable-order ODE solver (Matlab \texttt{ode113}) with a time step 100 times smaller than the discretization grid of the optimization problem. Observe that:
\begin{itemize}
\item Despite a very small tolerance and successful termination of the NLP solver, the collocation solution and interpolation of the solution exhibit large errors, leading to significant deviations to the state trajectories when the inputs are directly applied. In contrast, only minor discrepancies can be observed for the solutions represented using integrated residual minimization, on the same coarse grid with relatively low-order discretization. 
\item Although the constraints are implemented in the exact same way, the proposed method, to a greater extent, alleviates the issues of constraint violations inside the mesh intervals. This is because these constraint violations are often related to the large ODE defect errors  in-between  collocation points, which are directly dealt with by the residual minimization scheme.
\end{itemize}

\subsection{Aircraft Go-around in the Presence of Windshear}
\label{sec: AircraftExample}

\begin{figure}[tb]
\begin{center}
\includegraphics[width=0.9\columnwidth]{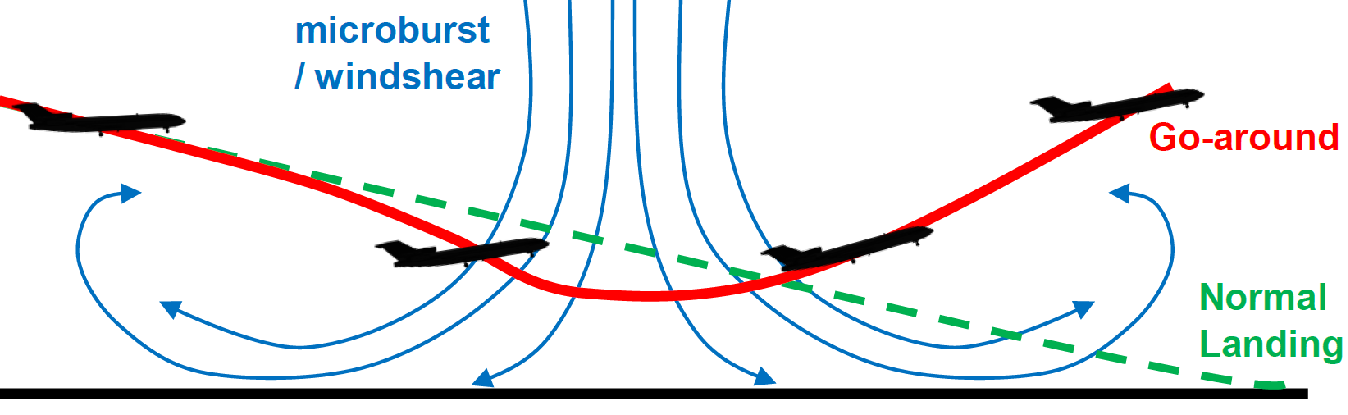}
\caption{Illustration of the aircraft go-around in windshear problem} 
\label{fig:AGIWSIllustration}
\end{center}
\end{figure}

Based on previous developments \cite{miele1988optimal}, a problem is presented in \cite{betts2010practical} where the aircraft needs to stay as high above the ground as possible after encountering a severe windshear during landing.  See the illustration in Figure~\ref{fig:AGIWSIllustration}. 

A simplified windshear model is used with wind speed contributions represented by a horizontal and a vertical component. Other details about the aerodynamic modelling, parameter values, simple bounds and boundary conditions are the same as in \cite{betts2010practical}. A static parameter $h_{min}$ is introduced to represent the minimum altitude. The objective is therefore to minimize $-h_{min}$ together with path constraint $h(t) \ge h_{min}$.

The angle of attack $\alpha$ is the actual control input to the physical system (aircraft); however, in order to implement a constraint on its rate of change, $\nu$ is introduced as angle of attack rate and serves as the control input with $\dot{\alpha}(t) =  \nu(t)$. This implementation is known to exhibit singular arc behaviour \cite{betts2010practical, nie2018should, rateconstraintArxiv}, leading to fluctuations and ringing phenomena in the solutions.

\begin{figure}[t]
\begin{center}
\includegraphics[width=\columnwidth]{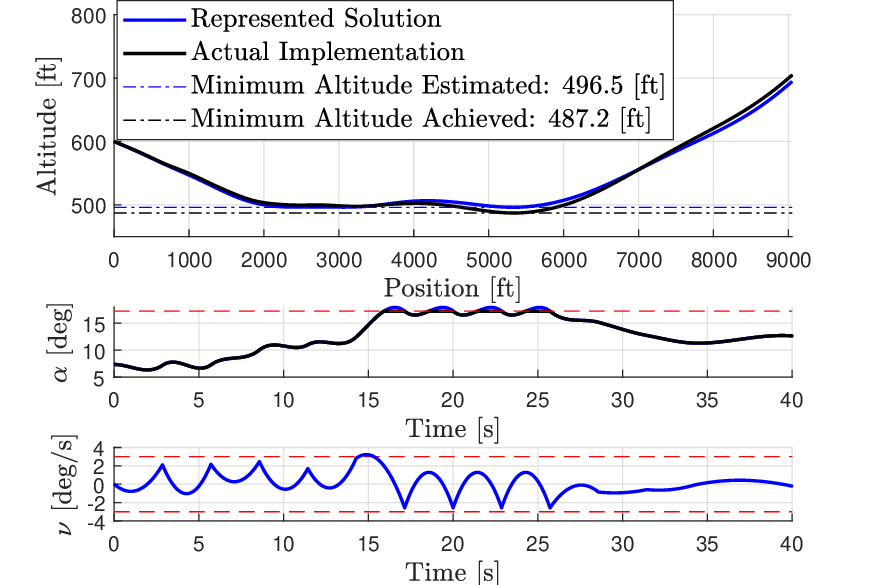}    
\caption{Solution to the aircraft go-around in the windshear problem, direct interpolation method for solution representation, direct collocation with Hermite-Simpson discretization, 15 mesh intervals} 
\label{fig:CollocationSolution}
\end{center}
\end{figure}

\begin{figure}[t]
\begin{center}
\includegraphics[width=\columnwidth]{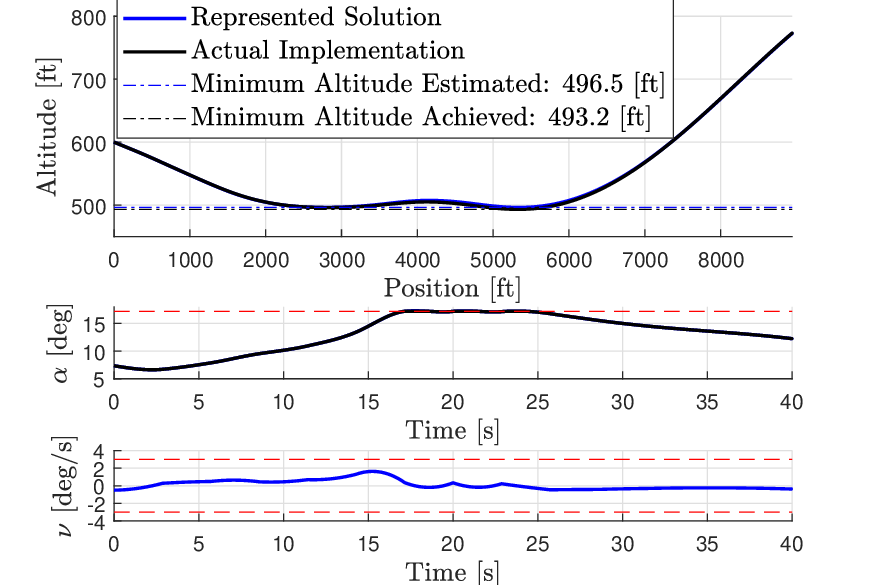}    
\caption{Solution to the aircraft go-around in the windshear problem, integrated residual minimization method for solution representation, direct collocation with Hermite-Simpson discretization, 15 mesh intervals} 
\label{fig:ResMinSolution}
\end{center}
\end{figure}

The solutions to this problem are collectively shown in Figures \ref{fig:CollocationSolution} and \ref{fig:ResMinSolution}. 
In addition to the advantages identified in the previous example, solution representation via integrated residual minimization have clear benefits in suppressing fluctuations. This ringing phenomenon is frequently observed in direct collocation solutions of singular control problems, as well as the directly interpolated solution trajectories.

\section{Conclusions}
\label{sec: conlclusions}

Although interpolation of collocation solutions can sometimes yield good results, it is very difficult to guarantee accuracy without posterior error assessments and mesh design iterations. As shown by the examples, despite successfully solving the NLP to negligibly small tolerance, the validity of the solution may still be questionable with large discrepancies. The flaws are rooted in collocation schemes, where the ODE defect errors are forced to zero at collocation points, regardless of the errors inside the intervals.  

The proposed solution representation method of integrated residual minimization fundamentally addresses this shortcoming by instead minimizing the integrated ODE residual error along the whole trajectory. As a result, solutions of higher accuracy are obtainable with the same discretization mesh, allowing the mesh to be relatively coarse. This benefit is clearly demonstrated with the example problems: despite being highly nonlinear, moderately complex and solved on a coarse low-order mesh, only minor differences are observed between the represented solution and actual implementation.


Since solving OCPs with numerical methods are essentially multi-objective optimization problems, one will inevitably face the trade-off between minimizing the objective (for optimality) and minimizing the discretization errors (for accuracy). By utilizing the objective value from a collocation solution, our proposed approach is yet to offer complete flexibility for the user in managing this trade-off process. Further development of this concept into a standalone method could potentially offer more benefits.

\bibliography{JournalBib2018} 
\bibliographystyle{ieeetr}

\end{document}